\documentclass{article}
\usepackage{amssymb} 
 \usepackage{amsmath}
\usepackage{amsthm}
\makeatletter
\def\@seccntformat#1{\csname named#1\endcsname\csname the#1\endcsname.\ }
\makeatother

\theoremstyle{theorem}

\theoremstyle{definition}

\def\ve#1{{\bf #1}}
\def\col{\mathop{\rm col}\nolimits}

\begin{document}

\centerline{\bf Two Eigenvectors for the price of one}

\rightline{Juan Tolosa}
\rightline{Stockton University}
\rightline{{\tt juan.tolosa@stockton.edu}}

\medskip

\begin{abstract}
Starting from a mistake done by a student, 
we discover an unexpected method of finding 
both eigenvectors for a $2\times2$ matrix with distinct eigenvalues in a single computation. 
We discuss a connection with the Cayley-Hamilton theorem, 
and show the corresponding generalization for a $3 \times 3$ matrix. 
The arguments should be understandable for strong linear algebra students. 
\end{abstract}

\section{Introduction. The $2\times2$ case.} 
It all started with a student's mistake. In a quiz, students were asked to find eigenvalues and eigenvectors
of the matrix
$$ A = 
\begin{pmatrix}
4&1\cr 2&5
\end{pmatrix}
$$
One student correctly found the two eigenvalues 
$\lambda_1 = 3$ and $\lambda_2=6$. She also correctly computed the matrix
$$ B_1 = A - \lambda_1 I =  \begin{pmatrix}
1&1\cr2&2\cr
\end{pmatrix}
$$
(where $I$ is the identity matrix), only instead of solving the homogeneous system
$$B_1 \ve v_1 = \ve 0 , $$
to find the corresponding eigenvector $\ve v_1 = {x \choose y}$ for $\lambda_1 = 3$, 
she picked a column of $B_1$ and (wrongly) declared that this was the required eigenvector, 
instead of picking, say, $\ve v_1 = {-1 \choose 1}$.

To my surprise, however, 
I realized that the chosen vector, $\ve v_2 = {1 \choose 2}$, is
actually an eigenvector for {\it the second eigenvalue\/} $\lambda_2 = 6$! 

At first I decided this might be just a coincidence. Upon checking the general case of a 
real-valued $2 \times 2$ matrix $A$ with distinct eigenvalues, this turned out to 
be true in general. The direct check could be a good exercise for a motivated student. 

We can, instead, seek a deeper understanding of the situation and, in doing so, 
discover a connection with an important result in linear algebra.
Namely, assuming $A$ has distinct eigenvalues $\lambda_1 \neq \lambda_2$, 
if we call 
$B_2 = A - \lambda_2 I$ the matrix we use to find the second eigenvector,
then the fact that the columns of $B_1$ turn out to be solutions $\ve v_2$ 
of $B_2 \ve v_2 = \ve 0$ is equivalent to saying that the matrix product $B_2 B_1$ is equal to 
the zero matrix. But if we think a little, we realize that 
$$ B_2 B_1 = (A - \lambda_2 I)(A - \lambda_1 I) = A^2 - (\lambda_1 + \lambda_2) A
+ \lambda_1 \lambda_2 I $$
is equal to 
$p(A)$, where 
$p(\lambda)$ is the characteristic polynomial of $A$, 
$p(\lambda) = \lambda^2 - ({\rm trace}\, A) \lambda + {\rm (determinant\,} A)$.
And the fact that $p(A)$ is indeed zero 
is the celebrated Cayley-Hamilton theorem, stating that every square matrix 
satisfies its own characteristic equation; 
see, for example, \cite{axler} or \cite{crilly}.

\smallskip

Thus, for a $2\times2$ matrix $A$ with distinct eigenvalues, 
the computations 
to find the first eigenvector $\lambda_1$  will also provide the second eigenvector for free! 
We just need to pick any {nonzero} column vector of $B_1 = A - \lambda_1 I$. 
Notice that, since we assume eigenvalues are distinct, then the kernel of 
$B_1$ (the corresponding eigenspace) has dimension 1, and the rank-nullity theorem 
guarantees the dimension of the column space of $B_1$ is also $1$, so $B_1$ 
is guaranteed to have at least one nonzero column.

\medskip

If $A$ has a double eigenvalue $\lambda_1 = \lambda_2$, 
then by Cayley-Hamilton $p(A) = B_1^2 = 0$, and we have two 
possible cases:
\begin{enumerate}
\item
 $B_1 = 0$ or, equivalently, $A$ is a multiple of the identity.
 Here all nonzero vectors are eigenvectors; 
 we can pick any two linearly independent vectors to get the basis of the eigenspace.
 \item
 $B_1 \neq 0$, in which case 
the eigenspace has dimension one;
any nonzero column of $B_1$ will provide an eigenvector. 
If we want a basis in which $A$ has canonical (Jordan) 
form, then, besides an eigenvector,  
we need to also pick a 
generalized eigenvector 
(see \cite{axler}, p.~66). 
Here is how we can get both.
Since $B_1^2 = 0$, one needs just pick any vector $\ve w$ 
not in the kernel of $B_1$
(that is, any nonzero vector which is not an eigenvector),
and compute $\ve v = B_1 \ve w$. Then $\ve v$ will be an eigenvector of $A$, 
with generalized eigenvector $\ve w$. Namely, we will have 
$$ A \ve v = \lambda \ve v \qquad\hbox{ and }\qquad
A \ve w = \lambda \ve w + \ve v . $$
For example, for the matrix
$$ A = \begin{pmatrix}
2&1\cr-1&4
\end{pmatrix} , $$
with double eigenvalue $\lambda = 3$, we have
$$ B = A - 3 I = \begin{pmatrix}
-1&1\cr-1&1
\end{pmatrix} . $$
if we pick $\ve w = (1, 0)$,
we get the eigenvector $ \ve v = B \ve w = (-1, -1)$, 
with generalized eigenvector $\ve w$.
\end{enumerate}

\medskip


\section{The $3 \times3$ case.} 
What about a $3 \times 3$ matrix $A$? 
First let us assume that $A$ has three distinct eigenvalues 
$\lambda_1, \lambda_2, \lambda_3$, and compute the corresponding matrices
$B_i = A - \lambda_i I$, $i=1,2,3$. 

Then the Cayley-Hamilton theorem will imply that 
$B_1 B_2 B_3 = 0$, so the nonzero columns 
of the product $B_2 B_3$ will be eigenvectors of 
$\lambda_1$. For example, the matrix 
$$ A =
\begin{pmatrix}
7&-4&-5\cr   3&-2&-3\cr  6&-4&-4
\end{pmatrix}
$$
has eigenvalues $\lambda_1 = 1$, $\lambda_2 = 2$, and $\lambda_3 = -2$. 
Here 
$$ B_2 = 
\begin{pmatrix}
5&-4&-5\cr  3&-4&-3\cr  6&-4&-6
\end{pmatrix} , \qquad
B_3 = 
\begin{pmatrix}
9&-4&-5\cr  3&0&-3\cr  6&-4&-2
\end{pmatrix} ,
$$
and
$$ B_2 B_3 = 
\begin{pmatrix}
3&0&-3 \cr  -3&0&3\cr   6&0&-6
\end{pmatrix} .
$$
The first and the last columns provide an eigenvector for $\lambda_1 = 1$, 
which can be simplified to $(1, -1, 2)$. 

Notice that $B_2 B_3$ must have at least one nonzero column. Indeed, if 
the three eigenvalues are distinct, then $p(\lambda)$ is the 
{\it minimal polynomial\/} of $A$ (\cite{axler}, p.~179), 
whereas $(\lambda - \lambda_2)(\lambda - \lambda_3)$ has smaller degree, so 
$B_2B_3 = (A - \lambda_2 I)(A- \lambda_3 I)$ cannot be the zero matrix. 

\smallskip

Assume now $A$ has a simple eigenvalue $\lambda_1$ and a double eigenvalue $\lambda_2$. 
Then we have two sub-cases. 

\begin{enumerate}

\item {\bf The eigenspace of $\lambda_2$ has dimension two.} 
Then the minimal polynomial of $A$ is $(\lambda - \lambda_1)(\lambda - \lambda_2)$, so that
not only $B_2^2 B_1 = 0$, but also $B_2 B_1 = 0$; 
therefore, the nonzero columns of $B_1$ are eigenvectors 
of the eigenvalue $\lambda_2$. Moreover, since $\lambda_1$ is a simple eigenvalue, the kernel of 
$B_1$ has dimension $1$; hence, the column space of $B_1$ has dimension $2$. In other words,
the column space of $B_1$ generates 
{\it the entire\/} (two-dimensional) eigenspace. 
So here you get {\it three\/} eigenvectors for the price of one!
As an example, consider the  matrix
$$ A =
\begin{pmatrix}
4&-9&-6\cr   -6&7&6\cr   12&-18&-14
\end{pmatrix}
$$
with $\lambda_1 = 1$ and 
$\lambda_2 = \lambda_3 = -2$. 
Here 
$$
B_1 =
\begin{pmatrix}
3&-9&-6\cr   -6&6&6\cr   12&-18&-15
\end{pmatrix}
$$
can be used to solve $B_1 \ve v = \ve 0$ and find the eigenvector
$\ve v_1 = (1,-1,2)$, corresponding to $\lambda_1 = 1$. 
Moreover,  
a direct check shows that {\it all\/} columns of $B_1$ are eigenvectors for 
$\lambda_2 = \lambda_3 = -2$.
Since $\dim(\col(B_1)) = 2$ (because $\lambda_1$ is a simple eigenvalue),
we conclude that the eigenspace of $\lambda_2$ is two-dimensional, and
the columns of $B_1$ generate this eigenspace. 
For example, the first two columns of $B$ (suitably divided by $3$) provide the 
missing eigenvectors 
$\ve v_2 = (1, -2, 4)$ and 
$\ve v_3 = (-3, 2, -6)$ for $\lambda_2$.

\item {\bf The eigenspace of $\lambda_2$ has dimension one.} 
Then the characteristic polynomial $p(\lambda) = (\lambda - \lambda_1)(\lambda - \lambda_2)^2$ 
is also the minimal polynomial. 
Hence, $B_2^2 B_1 = 0$ but $B_2 B_1 \neq 0$. 
As in the first case, the column space of $B_1$ has dimension 2. 
Moreover, in this case the kernel of $B_2$ has dimension 1, and the kernel of $B_2^2$ 
has dimension 2. 
(It has dimension at least 2, since $\col(B_1) \subset \ker(B_2^2)$, and cannot have dimension  
$3$, since otherwise $A$ would be a multiple of the identity, and $\lambda_2$ would be 
a triple eigenvalue.)
Consequently, the column space of $B_1$ generates the entire set of eigenvectors and 
generalized eigenvectors of $\lambda_2$. 
Hence, to find an eigenvector and generalized eigenvector for $\lambda_2$ 
we may proceed as follows:

\begin{itemize}
\item
Pick a nonzero column $\ve w$ of $B_1$ that is not an eigenvector of $B_2$, 
that is, such that $B_2 \ve w$ is nonzero. Such column must exist, since $\dim(\col(B_1))=2$
and $\dim(\ker(B_2)) = 1$.

\item
Compute $\ve v = B_2 \ve w$. 
Then $B_2 \ve v = B_2^2 \ve w = 0$, since $\ve w \in \col(B_1) = \ker (B_2^2)$.
We conclude that $\ve v$ is an eigenvector, and $\ve w$ is the 
corresponding generalized eigenvector, 
for $\lambda_2$.
\end{itemize}

As an example, consider
$$ A =
\begin{pmatrix}
5&-10&-7\cr
-6&7&6\cr
13 & -19 & -15
\end{pmatrix} ,
$$
with eigenvalues $\lambda_1 = 1$ and $\lambda_2 = \lambda_3 = -2$.
Here
$$ B_1 = A - \lambda_1 I =
\begin{pmatrix}
4 & -10 & -7\cr
-6&6&6 \cr
13&-19&-16
\end{pmatrix} .
$$
As usual, we use $B_1$ to find an eigenvector for $\lambda_1$, say,
$\ve v_1 = (1, -1, 2)$.

Next, 
a direct check shows that not all columns of $B_1$ are eigenvectors for $\lambda_2$ 
(in fact, {\it neither column\/} of $B_1$ is). 
This tells us that the eigenspace for $\lambda_2$ 
has dimension 1 (since otherwise {\it all\/} columns of $B_1$ would be eigenvectors).
To proceed, 
we can therefore pick any nonzero column of $B_1$, 
say, the first column:
$$ \ve w = 
\begin{pmatrix}
4\\ -6\\ 13
\end{pmatrix}, $$
and calculate
$$ \ve v_2 = B_2 \ve w = 
\begin{pmatrix}
-3 \\ 0\\ -3
\end{pmatrix} .
$$
We conclude that $\ve v_2$ is an eigenvector for $\lambda_2$, with corresponding 
generalized eigenvector $\ve w$. 
In the basis $\{ \ve v_1, \ve v_2, \ve w\}$, the matrix $A$ will have (Jordan) canonical form. 

\end{enumerate}

\medskip


\section{A triple eigenvalue.}
If $A$ is a $3\times3$-matrix with a triple eigenvalue $\lambda$, we have three possibilities 
(Let us again denote $B = A - \lambda I$).
\begin{enumerate}
\item
The eigenspace of $\lambda$ has dimension three, then $B = 0$, and all nonzero vectors are 
eigenvectors.
\item
The eigenspace of $\lambda$ has dimension one.
Then $B^2 \neq 0$ and has rank one; any nonzero column of $B^2$ will provide an eigenvector.
Moreover, if we want a basis in which $A$ will have canonical form, we can get two generalized vectors
as follows:
Pick any nonzero vector $\ve v_1$ that is not in the kernel of $B$; 
compute $\ve v_2 = B \ve v_1$, and $\ve v_3 = B \ve v_2$. 
Then $\ve v_3$ is an eigenvector of $A$, with two generalized eigenvectors $\ve v_2$ and $\ve v_1$.
\item
The eigenspace of $\lambda$ has dimension two. Then $B \neq 0$, $B^2 = 0$, and $B$ has rank one. 
Any nonzero column of $B$ will be an eigenvector of $A$. Moreover, if we choose any vector 
$\ve v$ not in the kernel of $B$, and compute $\ve w = B \ve v$, we will have both an eigenvector
$\ve w$ and a corresponding generalized vector $\ve v$.
Unfortunately, here we will still be missing a second (independent) eigenvector.
As it turns out, however, a careful analysis of the linear dependence of the columns of 
$B$ {\it does provide\/} two independent eigenvectors; see Section \ref{SpecialCase}.
\end{enumerate}

\medskip


\section{What if we don't want to solve any linear systems?}
\label{WhatIf}
It turns out in all cases 
we can find the eigenvectors without solving any linear systems!
Here is how we do it (the missing justifications are in previous sections).
\begin{itemize}
\item
{$2\times2$-matrix $A$}. 
\begin{enumerate}
\item
If the eigenvalues $\lambda_1, \lambda_2$ are distinct, compute $B_i = A - \lambda_i I$; 
any nonzero column of $B_1$ is an eigenvector for $\lambda_2$, and vice-versa.
\item If $\lambda_1 = \lambda_2$ and the eigenspace has dimension two, 
then $A$ is a multiple of the identity; 
pick any two linearly independent vectors as a basis for the eigenspace.
\item If $\lambda_1 = \lambda_2$ and the eigenspace has dimension one, 
then any nonzero column of $B_1$ will generate the eigenspace.
\end{enumerate}
\item
{$3\times 3$-matrix $A$}.
\begin{enumerate}
\item
If eigenvalues $\lambda_i$ are distinct, compute $B_i = A - \lambda_i I$, $i = 1, 2, 3$.
Then any nonzero column of $B_i B_j$ is an eigenvector for the third $\lambda_k$.
\item
If $\lambda_1 \neq \lambda_2 = \lambda_3$. 
Since $B_1 B_2^2 = 0$, then any nonzero column of $B_2^2$ spans the eigenspace for $\lambda_1$. 
Moreover,
\begin{itemize}
\item
If $B_2 B_1 = 0$, then the eigenspace of $\lambda_2$ has dimension two, and is spanned
by any two linearly independent columns of $B_1$.
\item
If $B_2 B_1 \neq 0$, then the eigenspace of $\lambda_2$ is one-dimensional, and is 
spanned by any nonzero column of $B_2 B_1$.
\end{itemize} 
\item
If $\lambda$ is a triple eigenvalue, we have three possibilities.
\begin{itemize}
\item
The eigenspace of $\lambda$ has dimension three. 
Then $A$ is a multiple of the identity, and any (nonzero) vector is an eigenvector.
\item
The eigenspace of $\lambda$ has dimension one.
Then $B^2 \neq 0$, and any nonzero column of $B^2$ will span this eigenspace.
\item
The eigenspace of $\lambda$ has dimension two.
Then $B \neq 0$, $B^2 = 0$, and $B$ has rank one. 
Any nonzero column $\ve v$ of $B$ will provide an eigenvector, but  
we are still missing a second eigenvector. 
It appears at first sight 
that the columns of $B$ provide no means 
of finding 
the entire eigenspace of $\lambda$ 
without solving a linear system. 
In the following Section we will show, however, that this is not so.

\end{itemize}
\end{enumerate}
\end{itemize}

\section{More on the special case}
\label{SpecialCase}
In Section \ref{WhatIf} we showed how to find the eigenspace of 
a $3\times3$ matrix $A$ 
using the columns of $B_i = A - \lambda_i I$ 
 in all cases but one, Case 3, 
when $A$ has a triple eigenvalue $\lambda$ with a two-dimensional eigenspace. 
Then the column space of 
$B = A - \lambda I$ is one dimensional, and provides one of the two independent
eigenvectors of the eigenspace of $\lambda$; 
the corresponding eigenvalue for $B$ is zero. 
The following result is straightforward.

{\narrower\narrower\it
A matrix $B$ has a triple eigenvalue zero, with two-dimensional eigenspace 
if, and only if, 
$B \neq 0$ and $B^2 = 0$. 
\par} 

For such a matrix the column space is one-dimensional, 
whereas the eigenspace is two-dimensional. 
Any nonzero column of $B$ provides an eigenvector, but we are missing a second one.

It turns out, however, that {\it the columns of $B$ 
do provide the two required eigenvectors\/}; it is just a matter of 
analyzing the specific linear dependence of the columns of $B$, as we will see below. 

The proofs of all the results stated in this section are skipped; 
they can be carried out by direct calculation.

\medskip

Since the column space of $B$ is one-dimensional, it follows that 
all columns of $B$ are proportional to one of them (any nonzero column).

However, this condition is not sufficient to guarantee that $B^2$ is the zero matrix,
as seen by the example
$$ B = \begin{pmatrix}
1&1&1\cr
1&1&1 \cr
1&1&1
\end{pmatrix} .
$$

To figure out sufficient conditions, and the structure of the eigenspace, 
let us divide the discussion into three cases.

\subsection
{Case 1. The first column of $B$ is nonzero.}

Then the matrix $B$ has the following structure, by columns:
\begin{equation}
\label{case1}
B = \left( \ve v \, \left\vert t \ve v \right \vert  s \ve v      \right)  , 
\end{equation} 
where 
$$ \ve v = 
\begin{pmatrix} x \cr y \cr  z
\end{pmatrix} $$ 
is the first column of $B$, which we assume to be nonzero, 
and $t, s$ are arbitrary real numbers.

Here we can prove the following fact.

{\narrower\narrower\it
The matrix $B$ satisfies $B^2 = 0$ if, and only if, 
$$ x + ty + sz = 0 . $$
Moreover, in this case, the eigenspace of the 
(triple) zero eigenvalue is spanned by the vectors
\begin{equation}
\label{case1vec}
 \ve v_1 = \begin{pmatrix} 
-t \cr 1 \cr 0
\end{pmatrix}
\qquad \hbox{ and } \qquad 
\ve v_2 = 
\begin{pmatrix}
-s \cr 0 \cr 1
\end{pmatrix}. 
\end{equation}
\par}

Thus, in Case 1, the matrix $B$ has the form (\ref{case1}), with
$$ \ve v = 
\begin{pmatrix} -ty - sz \cr y \cr  z
\end{pmatrix} , $$ 
and the eigenspace of zero is spanned by the vectors (\ref{case1vec}).

Such matrix $B$ constitutes a four-parameter family of matrices,
with arbitrary parameters $t, s, y, z$, with the only condition that not 
both $(y, z)$ are zero.

{\bf Example.}
An example with triple eigenvalue $3$: 
$$
A = \begin{pmatrix}
-2 & 5 & -10 \cr
-1 & 4 & -2 \cr
2  & -2 & 7
\end{pmatrix} ;
\qquad 
B = A - 3 I = \begin{pmatrix}
-5 & 5 & -10 \cr
-1 & 1 & -2 \cr
2  & -2 & 4
\end{pmatrix} ,
$$
so here $t = -1$ and $s = 2$. 
It follows that the eigenspace of the eigenvalue $3$ is spanned by
$
 \ve v_1 = \begin{pmatrix} 
1 \cr 1 \cr 0
\end{pmatrix}
$
and 
$
\ve v_2 = 
\begin{pmatrix}
-2 \cr 0 \cr 1
\end{pmatrix}. 
$
Of course, the columns of $B$ also belong to this eigenspace. 

Moreover, if needed, it is not hard to indicate a generalized eigenvector
for a general matrix in Case 1. Since the first column $\ve v$ is an eigenvector,
and since $B \ve e_1 = $ the first column of $B$, 
where $\ve e_1 = \begin{pmatrix} 1 \cr 0 \cr 0 \end{pmatrix}$, it follows that 
$\ve e_1$ is a generalized eigenvector, corresponding to the eigenvector $\ve v$.

\subsection
{Case 2. The first column of $B$ is zero, the second column is nonzero.}

Then the matrix $B$ has the following form, by columns:
\begin{equation}
\label{case2} 
B = \left( \ve0 \, \left\vert \ve v \right \vert  t \ve v      \right)  , 
\end{equation} 
where 
$$ \ve v = 
\begin{pmatrix} x \cr y \cr  z
\end{pmatrix} $$ 
is the second (nonzero) column of $B$, 
and $t$ is an arbitrary real number.

The following result holds in this case:

{\narrower\narrower\it
The matrix (\ref{case2})  satisfies $B^2 = 0$ if, and only if, 
$$y + tz = 0 . $$
Moreover, in this case, the eigenspace of the 
(triple) zero eigenvalue is spanned by the vectors
\begin{equation}
\label{case2vec}
 \ve v_1 = \begin{pmatrix} 
0 \cr -t \cr 1
\end{pmatrix}
\qquad \hbox{ and } \qquad 
\ve v_2 = \ve e_1 = 
\begin{pmatrix}
1 \cr 0 \cr 0
\end{pmatrix}. 
\end{equation}
\par}

It follows that in Case 2 the matrix $B$ has the form (\ref{case2}), with 
$$ \ve v = 
\begin{pmatrix} x \cr - tz \cr  z
\end{pmatrix} , $$ 
and the eigenspace of zero is spanned by the vectors (\ref{case2vec}).

Here we have a three-parameter family of matrices $B$, with parameters 
$t, x, z$, with the only condition that $(x, z)$ are not simultaneously zero. 

{\bf Example.} A matrix with triple eigenvalue $2$: 
$$
A = \begin{pmatrix}
2 & -1 & 2 \cr
0 & -4 & 12 \cr
0 & -3 & 8
\end{pmatrix}
; \qquad
B = A - 2 I = \begin{pmatrix}
0 & -1 & 2 \cr
0 & -6 &  12 \cr
0 & -3 & 6
\end{pmatrix} ,
$$
so $t = -2$, and
therefore the eigenspace of $\lambda = 2$ is spanned by
$
 \ve v_1 = \begin{pmatrix} 
0 \cr 2 \cr 1
\end{pmatrix}
$ 
and 
$ 
\ve e_1 = 
\begin{pmatrix}
1 \cr 0 \cr 0
\end{pmatrix}. 
$
Moreover, 
$\ve e_2 = \begin{pmatrix} 0 \cr 1 \cr 0 \end{pmatrix}$
is a generalized vector corresponding to the eigenvector
$\ve v = \begin{pmatrix} -1 \cr -6 \cr -3 \end{pmatrix}$,
since $B \ve e_2 = \ve v  = $ second column of $B$.

\subsection
{Case 3. The first two columns of $B$ are zero, the last column is nonzero.}

Here the following holds.

{\narrower\narrower\it
In Case 3, 
the matrix $B$ satisfies $B^2 = 0$ if, and only if, it has the form
$$B = \begin{pmatrix}
0 & 0 & x \cr
0 & 0 & y \cr
0 & 0 & 0 
\end{pmatrix} ,
$$
where $x, y$ are arbitrary numbers, not both zero. 
Moreover, in this case, the eigenspace of the 
(triple) zero eigenvalue is spanned by the vectors
\begin{equation}
 \ve e_1 = \begin{pmatrix} 
1 \cr 0 \cr 0
\end{pmatrix}
\qquad \hbox{ and } \qquad 
\ve e_2 =  
\begin{pmatrix}
0 \cr 1 \cr 0
\end{pmatrix}. 
\end{equation}
\par}

Finally,
$ \ve e_3 = \begin{pmatrix} 
0 \cr 0 \cr 1
\end{pmatrix}$
is a generalized vector for the eigenvector 
$ \ve v =\begin{pmatrix} 
x \cr y \cr 0
\end{pmatrix} . $

%
%

\section{Acknowledgment.}
I am grateful to the editors of the {\it College Mathematics Journal\/} 
for most of the content of Section \ref{WhatIf}.


%


\vfill\eject


\begin{thebibliography}{1}
\bibitem{axler} Axler, S. (1997). \textit{Linear Algebra Done Right.} New York: Springer; 
p.~173 and p.~207.

\bibitem{crilly} Crilly, T. (1992). A Gemstone in Matrix Algebra \textit{Math. Gazette} 76(475): 182--188. doi.org/10.2307/3620391.
\end{thebibliography}
\end{document}